\documentclass[11pt,reqno]{amsart}

\usepackage{cite,a4wide,amsmath,amssymb,graphicx,enumitem,amsaddr,tikz,pgfplots,hhline,epstopdf,multirow,multicol,soul}
\usepackage[hang,flushmargin]{footmisc}

\usetikzlibrary{decorations.markings,patterns}

\newtheorem{theorem}{Theorem}

\setlist[enumerate]{leftmargin=8mm}
\setlist[itemize]{leftmargin=8mm}

\newcommand{\cR}{\mathcal{R}}

\newcommand{\e}{\text{e}}

\newcommand{\ds}{\frac{\text{d}S}{\text{d}t}}
\newcommand{\di}{\frac{\text{d}I}{\text{d}t}}
\newcommand{\dr}{\frac{\text{d}R}{\text{d}t}}
\newcommand{\de}{\frac{\text{d}E}{\text{d}t}}
\newcommand{\dv}{\frac{\text{d}V}{\text{d}t}}
\newcommand{\dq}{\frac{\text{d}Q}{\text{d}t}}
\newcommand{\dha}{\frac{\text{d}H}{\text{d}t}}

\DeclareMathOperator{\sgn}{sgn}

\definecolor{darkgreen}{RGB}{34,177,76}

\begin{document}
\title{From pandemic to a new normal: strategies to optimise governmental interventions in Indonesia based on an SVEIQHR-type mathematical model}

\author{\small Benny Yong, Jonathan Hoseana, Livia Owen}
\address{\normalfont\small Department of Mathematics, Parahyangan Catholic University, Bandung 40141, Indonesia}
\email{benny\_y@unpar.ac.id\textnormal{, }j.hoseana@unpar.ac.id\textnormal{, }livia.owen@unpar.ac.id}
\date{}

\begin{abstract}
There are five different forms of intervention presently realised by the Indonesian government in an effort to end the COVID-19 pandemic: vaccinations, social restrictions, tracings, testings, and treatments. In this paper, we construct an SVEIQHR-type mathematical model for the disease's spread in the country, which incorporates as parameters the rates of the above interventions, as well as the vaccine's efficacy. We determine the model's equilibria and basic reproduction number. Using the model, we formulate strategies by which the interventions should be realised in order to optimise their impact. The results show that, in a disease-free state, when the number of new cases rises, the best strategy is to implement social restrictions, whereas in an endemic state, if a near-lockdown policy is undesirable, carrying out vaccinations is the best strategy; however, efforts should be aimed not primarily towards increasing the vaccination rate, but towards the use of high-efficacy vaccines.



\end{abstract}

\maketitle

\section{Introduction}\label{section:Introduction}

Declared to be a pandemic by WHO on 11 March 2020 \cite[page 2]{Saxena}, the \textit{coronavirus disease 2019} (COVID-19), reportedly originating from a seafood market in China \cite[page 13--14]{Saxena}, has continued to be a global concern, with over 392 million cases recorded worldwide as of 6 February 2022 \cite{WHO}. In most countries, the initial evolution of the daily number of new cases is characterised by several successive waves \cite{MunosFernandezSeoaneSeoaneSepulveda}, the end of each wave seemingly indicating successfulness of certain eradicative interventions. Such successfulness was however only temporal in many countries where there occurred subsequent ---often larger--- waves, of which the emergence of new variants proved to be a principal cause \cite{Page}.

In Indonesia, since mid 2020, amid various disruptions caused by the disease, the government has popularised the term ``\textit{a new normal}'' to refer to a desired form of post-pandemic life \cite{Adjie}. Its realisability, however, remained unclear. Indeed, the largest pre-omicron wave unfolded in mid 2021, the number of new cases reaching a new maximum of 56,767 on 15 July 2021 \cite{AriefAdjiSuharto}. The omicron variant then entered the country in late November 2021 \cite{RanggasariBhwana}, before a subsequent wave emerged in January 2022 and the aforementioned maximum was surpassed as early as 15 February 2022, with 57,049 new cases \cite{Rayda}. Keeping the aim towards a new normal, the government of the country has been realising concrete eradicative interventions in the following five different forms.

\begin{enumerate}[itemsep=3pt]
\item \textit{Vaccinations}. Indonesia's national vaccination programme commenced on 13 January 2021 \cite{Triwardani}, the first vaccinated citizen being the president, Joko Widodo, who received on the day a shot of Coronavac, a vaccine developed by China's Sinovac Biotech, approved for emergency use by the country's Food and Drug Monitoring Agency (BPOM) only two days earlier \cite{Syakriah}. The programme's progress has been tangible: by early 2022, around 45\% and 21\% of the population have been fully and partially vaccinated, respectively \cite{WHOI19Jan2022}, and 11 different vaccines with varying levels of efficacy were granted approval \cite{COVID19VaccineTracker}. The reception of booster shots has also been urged \cite{Afifa}. Vaccinated citizens are provided with waivers from a number of health-related requirements\footnote{There is a true risk of such a policy, which will appear in a later discussion.} \cite{JunidaIhsanSuharto}.
\item \textit{Social restrictions}. Besides quarantine regulations for citizens returning from abroad \cite{IndonesiaExpat2}, the government, in an effort to control the disease's transmission level, has set out four different levels of large-scale social restrictions (PPKM) \cite{DepartemenDalamNegeri1,DepartemenDalamNegeri2,SaptoyoKurniawan,DoubleM}. Each level of social restrictions (1 to 4) defines a specific degree to which schools, shopping centres, public transport, etc., may operate. At any given time, every region is to implement one of these four restriction levels, carefully determined by the authority based on the region's present situation, using several indicators such as transmission and vaccination levels \cite{MenteriKesehatan}. In particular, responding to the omicron wave, the government has re-raised the levels of social restrictions in various regions including Greater Jakarta, from 2 to 3 \cite{Sucipto}.
\item \textit{Tracings}. Contact tracing, a procedure of interviewing a newly-diagnosed patient with the aim of identifying people who have been in contact with the patient within the last few days \cite{WHO2}, has also been in operation, albeit initially at a suboptimal level. (In fact, one of the reasons for the re-raising of the social restrictions level at the start of the omicron wave was the lack of tracings \cite{Sucipto}.) Recently, however, there has been some increase in awareness and government's effort towards contact tracings, as means to break transmission chains \cite{Hutasoit}.
\item \textit{Testings}. Since many COVID-19 patients are asymptomatic, determining whether a person ---including travellers and those who have been in recent contact with a patient--- is infected is best done via medical testings. At least three types of tests for COVID-19 are available in the country: the polymerase chain reaction test (PCR), the rapid antibody test, and the rapid antigen test \cite{NewsDesk}.
\item \textit{Treatments}. Being a populous country, Indonesia has faced some considerable difficulty in optimising medical treatments for COVID-19 patients as hospitals became increasingly crowded during the country's largest pre-omicron wave in mid 2021 \cite{WHOI23Jun2021}. A number of makeshift hospitals were set up, so as to keep the overall bed-occupancy rate ---and thus the quality of treatment--- within a safe level \cite{WHOI22Jul2020}. At the start of the omicron wave, the government has, in addition, looked into providing citizens with antiviral medicines, securing 400,000 tablets of molnupiravir \cite{Shofa}.
\end{enumerate}

In \cite{YongOwenHoseana}, we have constructed a SIR-type mathematical model for the spread of COVID-19 in Indonesia, which incorporates ---among others--- a parameter measuring the aforementioned hospitals' bed-occupancy rate. We have also used this model to design a quantitative method for determining the appropriate level(s) of social restrictions to be enforced in Jakarta at any given time \cite{YongHoseanaOwen}. In this paper, taking into consideration the above five forms of intervention, as well as the idea proposed in \cite[section 4]{YongOwenHoseana} of incorporating more compartments and the possibility of reinfection, we aim to construct a new model which is more comprehensive and realistic, with the hope of formulating strategies by which the above forms of intervention should be realised in order to optimise their impact, so that a new normal can be embraced as soon as possible.

The incorporation of additional compartments implies that the present model is no longer SIR-type. Indeed, we shall take into account, at any given time $t\geqslant 0$, the numbers $S=S(t)$ of (unvaccinated) \textit{susceptible} individuals, $V=V(t)$ of (susceptible) \textit{vaccinated} individuals, $E=E(t)$ of (non-transmitting) \textit{exposed} individuals, $I=I(t)$ of \textit{infected} individuals, $Q=Q(t)$ of \textit{quarantined} individuals, $H=H(t)$ of \textit{hospitalised} individuals, and $R=R(t)$ of \textit{recovered} individuals, thereby building a seven-compartment SVEIQHR-type model. We assume that all quarantines are centralised, so that the quarantined ---as well as the hospitalised--- individuals are never in contact with the susceptible and vaccinated individuals, meaning that only the infected individuals transmit the disease. We also assume that social restrictions are waived for vaccinated citizens\footnote{This could be an imprudent policy; see subsection \ref{subsec:R0}.}; see, e.g., \cite{JunidaIhsanSuharto}. The above five forms of intervention shall be incorporated to the model as parameters $u_1$, $u_2$, $u_3$, $u_4$, and $u_5$, all belonging to $[0,1]$, which represent, respectively, the rates of vaccine, mobility, contact-tracing, rapid-testing, and treatment interventions\footnote{Thus, $u_2=0$ represents normal mobility, while $u_2=1$ represents a total lockdown. For a method to estimate the value of $u_2$ representing each of the four levels of social restrictions, see subsection \ref{subsec:R0}.}. We shall also incorporate a parameter $\delta\in[0,1]$ representing the vaccine efficacy. These and all other parameters, together with their values used in our numerical analysis, are described in Table \ref{tab:parameters}.

\begin{figure}
\centering
\scalebox{0.85}{\begin{tikzpicture}
\node (muI) at (6.65,-1.75) {$\mu I$};
\node (mu'I) at (7.35,-1.75) {$\mu' I$};
\draw[<-,thick] (muI) edge (6.65,0);
\draw[<-,thick] (mu'I) edge (7.35,0);

\node (muQ) at (7,-5.25) {$\mu Q$};
\draw[<-,thick] (muQ) edge (7,-3.5);

\node (muH) at (10.2,-5.25) {$\mu H\,$};
\node (mu'H) at (10.8,-5.25) {$\,\mu' H$};
\draw[<-,thick] (muH) edge (10.2,-3.5);
\draw[<-,thick] (mu'H) edge (10.8,-3.5);

\node (muS) at (0.2,-1.75) {$\mu S$};
\draw[<-,thick] (muS) edge (0.2,0);

\node (muR) at (10.7,-1.75) {$\mu R$};
\draw[<-,thick] (muR) edge (10.7,0);

\draw[dashed,thick] (7,0) -- (7,1) -- (0.2,1) -- (0.2,0.4875);

\node[rectangle,fill=white,draw=black,thick,minimum size=0.95cm] (S) at (0,0) {S};
\node[rectangle,fill=white,draw=black,thick,minimum size=0.95cm] (E) at (3.5,0) {E};
\node[rectangle,fill=white,draw=black,thick,minimum size=0.95cm] (I) at (7,0) {I};
\node[rectangle,fill=white,draw=black,thick,minimum size=0.95cm] (R) at (10.5,0) {R};
\node[rectangle,fill=white,draw=black,thick,minimum size=0.95cm] (V) at (0,-3.5) {V};
\node[rectangle,fill=white,draw=black,thick,minimum size=0.95cm] (Q) at (7,-3.5) {Q};
\node[rectangle,fill=white,draw=black,thick,minimum size=0.95cm] (H) at (10.5,-3.5) {H};
\draw[dashed,thick,->] (3.5,1) -- (E);
\node (lambda) at (-1.75,0) {$\lambda$};
\draw[->,thick] (lambda) edge (S);
\node (muE) at (3.5,-1.75) {$\mu E$};
\draw[<-,thick] (muE) edge (E);
\node (muV) at (0,-5.25) {$\mu V$};
\draw[<-,thick] (muV) edge (V);
\draw[->,thick] (S) edge node[above] {$\displaystyle\left(1-u_2\right)\beta S I$} (E);
\draw[->,thick] (E) edge node[above] {$\theta E$} (I);
\draw[->,thick] (I) edge node[above] {$\gamma I$} (R);
\draw[->,thick] (R) -- (10.5,1.5) -- (-0.2,1.5) -- (-0.2,0.4875);
\node[above] at (5.25,1.5) {$\alpha R$};
\draw[->,thick] (S) edge [transform canvas={xshift=-2mm}] node[left] {$u_1 S$} (V);
\draw[->,thick] (V) edge node[right,pos=0.25] {$(1-\delta)\beta V I$} (E);
\draw[->,thick] (E) edge node[right] {$u_3E$} (Q);
\draw[->,thick] (I) edge node[right,pos=0.7] {$u_4I$} (Q);
\draw[->,thick] (Q) edge node[right,pos=0.2] {$\kappa Q$} (R);
\draw[->,thick] (I) edge node[right,pos=0.2] {$u_5 I$} (H);
\draw[->,thick] (H) edge[transform canvas={xshift=-2mm}] node[left] {$\varphi H$} (R);
\draw[->,thick] (Q) edge node[above] {$\tau Q$} (H);
\node (lambda') at (5.25,-3.5) {$\lambda'$};
\draw[->,thick] (lambda') edge (Q);
\end{tikzpicture}}
\caption{\label{fig:presentmodel} The compartment diagram of our SVEIQHR-type model.}
\end{figure}
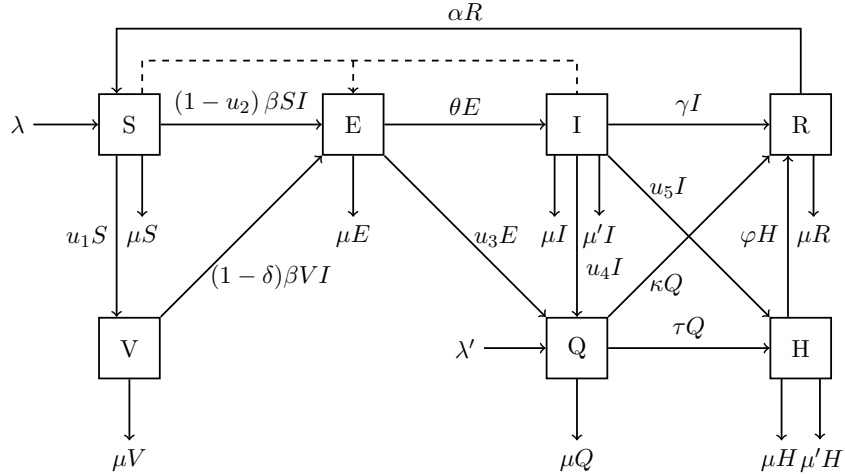

Let us now construct the model itself, by detailing the changes assumed to be experienced at any given time by each of the above seven time-dependent variables, which are summarised in the compartment diagram in Figure \ref{fig:presentmodel}.
\begin{enumerate}[itemsep=3pt]
\item The number $S$ of susceptible individuals increases due to the entry of newborns at the rate $\lambda>0$ and of recovered individuals at the rate $\alpha R$, where $\alpha>0$, and decreases due to the exit of those who become exposed at the rate $\left(1-u_2\right)\beta SI$, where $u_2\in[0,1]$ and $\beta>0$, vaccinated at the rate $u_1S$, where $u_1\in[0,1]$, and dead at the rate $\mu S$, where $\mu>0$.
\item The number $V$ of vaccinated individuals increases due to the entry of susceptible individuals at the rate $u_1S$, and decreases due to the exit of those who become exposed at the rate $(1-\delta)\beta VI$, where $\delta\in[0,1]$, and dead at the rate $\mu V$.
\item The number $E$ of exposed individuals increase due to the entry of susceptible individuals at the rate $\left(1-u_2\right)\beta SI$ and of vaccinated individuals at the rate $(1-\delta)\beta VI$, and decreases due to the exit of those who become infected at the rate $\theta E$, where $\theta>0$, quarantined at the rate $u_3 E$, where $u_3\in[0,1]$, and dead at the rate $\mu E$.
\item The number $I$ of infected individuals increase due to the entry of exposed individuals at the rate $\theta E$, and decreases due to the exit of those who become recovered at the rate $\gamma I$, where $\gamma>0$, quarantined at the rate $u_4 I$, where $u_4\in[0,1]$, hospitalised at the rate $u_5 I$, where $u_5\in[0,1]$, and dead at the rate $\left(\mu+\mu'\right)I$, where $\mu'>0$.
\item The number $Q$ of quarantined individuals increase due to the entry of foreigners at the rate $\lambda'>0$, of exposed individuals at the rate $u_3E$, and of infected individuals at the rate $u_4I$, and decreases due to the exit of those who become recovered at the rate $\kappa Q$, where $\kappa>0$, hospitalised at the rate $\tau Q$, where $\tau>0$, and dead at the rate $\mu Q$.
\item The number $H$ of hospitalised individuals increase due to the entry of quarantined individuals at the rate $\kappa Q$ and of infected individuals at the rate $u_5 I$, and decreases due to the exit of those who become recovered at the rate $\varphi H$, where $\varphi>0$, and dead at the rate $\left(\mu+\mu'\right)H$.
\item The number $R$ of recovered individuals increases due to the entry of infected individuals at the rate $\gamma I$, quarantined individuals at the rate $\kappa Q$, and hospitalised individuals at the rate $\varphi H$, and decreases due to the exit of those who become susceptible at the rate $\alpha R$ and dead at the rate $\mu R$.
\end{enumerate}

\begin{table}\centering
	\renewcommand{\arraystretch}{1.9}\renewcommand{\tabcolsep}{2pt}
		\scalebox{0.8}{\begin{tabular}{|c|c|c|c|c|} 
			\hline
			Parameter & Description & Unit & 
\begin{minipage}{2cm}\centering
Value for simulation
\end{minipage} & Source\\
			\hhline{|=|=|=|=|=|}
			$\lambda$ & recruitment rate of newborns & individual/day & $\displaystyle\frac{273523621}{65 \times 365}$ & \begin{minipage}{2.5cm}\centering
estimated as $\mu\, N(0)$ \cite{WorldBank}
\end{minipage}\\\hline
			$\lambda'$ & recruitment rate of foreigners & individual/day & 3000 & \cite{IndonesiaExpat}\\\hline
			$\mu$ & natural death rate & 1/day & $\displaystyle\frac{1}{65 \times 365}$ & \cite{AldilaEtAl1}\\\hline
			$\mu'$ & death rate increment due to COVID-19 & 1/day & 0.0291 & \cite{GugusTugas}\\\hline
			$\beta$ & transmission coefficient & 1/(individual\,$\times$\,day) & $4.74396\times 10^{-8}$ & \cite{AldilaEtAl1}\\\hline
			$\delta$ & vaccine efficacy & dimensionless & \multicolumn{2}{c|}{see subsection \ref{subsec:R0}}\\\hline
			$\alpha$ & temporary immunity rate & 1/day & 0.011 & \cite{ShakhanySalimifard}\\\hline
			$\theta$ & incubation rate & 1/day & 0.4 & \cite{PremEtAl}\\\hline
		    $\gamma$ & recovery rate of infected individuals & 1/day & 0.1 & \cite{FergusonEtAl}\\\hline
			$\varphi$ & recovery rate of hospitalised individuals & 1/day & 0.8198 & \cite{GugusTugas}\\\hline
			$\kappa$ & recovery rate of quarantined individuals & 1/day & 0.1 & \cite{FergusonEtAl}\\\hline
			$\tau$ & hospitalisation rate of quarantined individuals & 1/day & 0.01 & \cite{AldilaEtAl2}\\\hline
			$u_1$ & vaccination rate & 1/day & 0.4 & \cite{DiagneEtAl}\\\hline
			$u_2$ & mobility intervention rate & dimensionless & \multicolumn{2}{c|}{see subsection \ref{subsec:R0}}\\\hline
			$u_3$ & contact-tracing intervention rate & 1/day & 0.5 & assumed\\\hline
			$u_4$ & rapid-testing intervention rate & 1/day & 0.3 & assumed\\\hline
			$u_5$ & treatment intervention rate & 1/day & 0.0833 & \cite{BabaeiEtAl}\\\hline
		\end{tabular}}\smallskip
		\caption{\label{tab:parameters}Parameters used in the model \eqref{eq:model} and their values chosen for our numerical simulations (section \ref{section:analysis}).}
\end{table}

\noindent We therefore obtain the model
\begin{equation}\label{eq:model}
\left\{\begin{array}{rcl}
\displaystyle\ds &=& \displaystyle\lambda +\alpha R - \left(1-u_2\right)\beta SI -u_1S-\mu S,\\[0.3cm]
\displaystyle\dv &=& \displaystyle u_1 S-(1-\delta)\beta VI-\mu V,\\[0.3cm]
\displaystyle\de &=& \displaystyle\left(1-u_2\right)\beta SI  - \theta E+(1-\delta)\beta VI-u_3E-\mu E,\\[0.3cm]
\displaystyle\di &=& \displaystyle\theta E-\gamma I-u_4 I-u_5 I -\mu I -\mu' I,\\[0.3cm]
\displaystyle\dq &=& \displaystyle\lambda'+u_3E+u_4I-\kappa Q-\tau Q-\mu Q,\\[0.3cm]
\displaystyle\dha &=& \displaystyle\tau Q+u_5 I-\varphi H-\mu H-\mu'H,\\[0.3cm]
\displaystyle\dr &=& \displaystyle\gamma I-\alpha R+\kappa Q+\varphi H-\mu R.
\end{array}\right.
\end{equation}
 
The rest of the paper is organised as follows. In the upcoming section \ref{sec:equilibria}, we analyse the model \eqref{eq:model} dynamically. We first establish the non-negativity and boundedness of its solutions, and determine a subdomain which is positively invariant under the model (subsection \ref{subsec:boundedness}). Next, we show that, for every set of parameter values, the model possesses a unique disease-free equilibrium, and determine an explicit expression of this equilibrium (subsection \ref{subsec:DFER0}). We also derive the model's basic reproduction number $\cR_0$ and show that, if $\cR_0<1$, the disease-free equilibrium is stable, whereas if $\cR_0>1$, the disease-free equilibrium is unstable and a unique positive endemic equilibrium exists (subsections \ref{subsec:DFER0} and \ref{subsec:EE}).

As the algebraic computations required to establish further dynamical properties of the model ---such as the endemic equilibrium's stability--- appears to be inaccessibly complicated, we shift from analytical to numerical methods (section \ref{section:analysis}), whose flexibility allows us to achieve our ultimate goal: formulating strategies by which the aforementioned forms of governmental intervention (vaccinations, social restrictions, tracings, testings, and treatments) should be implemented for an optimal impact. The first stage of our analysis yields results which strongly point towards vaccinations, and more specifically, towards the importance of a high vaccine efficacy, in addition to the necessity of unwaiving social restrictions for vaccinated citizens (subsection \ref{subsec:R0}). This is confirmed quantitatively in our second stage (subsection \ref{subsec:sensitivity}) via sensitivity analysis, from which we conclude that the optimal intervention strategy is to implement social restrictions in the case of $\cR_0<1$, and, if a lockdown is undesirable, vaccinations using high-efficacy vaccines in the case of $\cR_0>1$. These conclusions are reasserted in section \ref{section:conclusions}, where we also describe a number of ways in which the model \eqref{eq:model} could be modified for further research.\smallskip

\section{Dynamical analysis}\label{sec:equilibria}

Let us first analyse the model \eqref{eq:model} from the viewpoint of dynamical systems theory; see \cite{Robinson,Martcheva} for background. First, we establish the non-negativity and boundedness of the model's solutions associated to non-negative initial conditions, and the positive-invariance of a bounded subdomain (Theorem \ref{thm:boundedness}). Subsequently, we show that the model has a unique disease-free equilibrium for every set of parameter values, which is stable if $\cR_0<1$ and unstable if $\cR_0>1$, where $\cR_0$ is the model's basic reproduction number (Theorem \ref{thm:DFE}). Finally, we show that in the case of $\cR_0>1$, in which the model's solutions do not approach the disease-free equilibrium, a unique positive endemic equilibrium exists (Theorem \ref{thm:EE}).\smallskip

\subsection{Non-negativity and boundedness of solutions}\label{subsec:boundedness} Let us first establish the non-negativity and boundedness of the solutions of the model \eqref{eq:model} associated to non-negative initial conditions. Let $$\left(S(0),V(0),E(0),I(0),Q(0),H(0),R(0)\right)\in[0,\infty)^7$$ be such an initial condition, and let $\left(S(t),V(t),E(t),I(t),Q(t),H(t),R(t)\right)$ be the solution associated to this initial condition. For every $t^\ast\geqslant 0$ satisfying $S\left(t^\ast\right)=0$, we have, from the model's first equation,
$$\left.\ds\right|_{t=t^\ast} = \lambda + \alpha R\left(t^\ast\right) > 0,$$
which means that the function $S$ is increasing at $t^\ast$. Since $S(0)\geqslant 0$, it follows that $S(t)\geqslant 0$ for every $t\geqslant 0$. Similar arguments show that
$$V(t)\geqslant0,\quad E(t)\geqslant0,\quad I(t)\geqslant0,\quad Q(t)\geqslant0,\quad H(t)\geqslant0,\quad \text{and}\quad R(t)\geqslant0$$
for every $t\geqslant 0$.

Next, adding all equations in \eqref{eq:model}, one obtains that the time-dependent total population $N:=S+V+E+I+Q+H+R$ satisfies
$$\frac{\text{d}N(t)}{\text{d}t}\leqslant\lambda +\lambda'-\mu N(t),$$
which is equivalent to
\begin{equation}\label{eq:slope}
\frac{\text{d}}{\text{d}t}\left(N(t)\e^{\mu t}\right)\leqslant\frac{\text{d}}{\text{d}t}\left(\frac{\lambda+\lambda'}{\mu}\e^{\mu t}+N(0)-\frac{\lambda+\lambda'}{\mu}\right).
\end{equation}
Now, the functions $N(t)\e^{\mu t}$ and $\left(\left(\lambda+\lambda'\right)/\mu\right)\e^{\mu t}+N(0)-\left(\lambda+\lambda'\right)/\mu$ have the same value at $t=0$, namely, $N(0)$, and, by \eqref{eq:slope}, at every point, the slope of the former function does not exceed that of the latter function. Consequently, for every $t\geqslant 0$ we have
$$N(t)\e^{\mu t}\leqslant \frac{\lambda+\lambda'}{\mu}\e^{\mu t}+N(0)-\frac{\lambda+\lambda'}{\mu},$$
i.e.,
$$N(t)\leqslant\frac{\lambda+\lambda'}{\mu}+\left(N(0)-\frac{\lambda+\lambda'}{\mu}\right)\e^{-\mu t}.$$

This implies that $$\lim_{t\to\infty} N(t)\leqslant \frac{\lambda+\lambda'}{\mu},$$
i.e., that the solution $\left(S(t),V(t),E(t),I(t),Q(t),H(t),R(t)\right)$ is bounded, and that the subset
$$\mathcal{D}:=\left\{(S,V,E,I,Q,H,R)\in[0,\infty)^7:S+V+E+I+Q+H+R\leqslant\frac{\lambda+\lambda'}{\mu}\right\}\subseteq [0,\infty)^7$$
\sloppy is positively invariant \cite[Definition 4.4]{Robinson} under the model.

We summarise our results in the following theorem.\smallskip

\begin{theorem}\label{thm:boundedness}\
\begin{enumerate}
\item[\textnormal{(1)}] Every solution of the model \eqref{eq:model} associated to an initial condition in $[0,\infty)^7$ is bounded and remains forever in $[0,\infty)^7$. 
\item[\textnormal{(2)}] Every solution of the model \eqref{eq:model} associated to an initial condition in $\mathcal{D}$ remains forever in $\mathcal{D}$.
\end{enumerate}
\end{theorem}\smallskip


\subsection{Disease-free equilibrium and basic reproduction number}\label{subsec:DFER0}

Let us now study the equilibria of the model \eqref{eq:model}, i.e., the solutions of the system
\begin{equation}\label{eq:equilibria}
\left\{\begin{array}{rcl}
\displaystyle\lambda +\alpha R - \left(1-u_2\right)\beta SI -u_1S-\mu S&=&0,\\[0.15cm]
\displaystyle u_1 S-(1-\delta)\beta VI-\mu V&=&0,\\[0.15cm]
\displaystyle\left(1-u_2\right)\beta SI  - \theta E+(1-\delta)\beta VI-u_3E-\mu E&=&0,\\[0.15cm]
\displaystyle\theta E-\gamma I-u_4 I-u_5 I -\mu I -\mu' I&=&0,\\[0.15cm]
\displaystyle\lambda'+u_3E+u_4I-\kappa Q-\tau Q-\mu Q&=&0,\\[0.15cm]
\displaystyle\tau Q+u_5 I-\varphi H-\mu H-\mu'H&=&0,\\[0.15cm]
\displaystyle\gamma I-\alpha R+\kappa Q+\varphi H-\mu R&=&0.
\end{array}\right.
\end{equation}
We shall begin by showing that, for every set of parameter values, the model possesses a unique disease-free equilibrium, which admits an explicit description, and relating its stability to the model's basic reproduction number.

Let $\mathbf{e}_0=\left(S_0,V_0,E_0,I_0,Q_0,H_0,R_0\right)$ be a disease-free equilibrium of the model \eqref{eq:model}, i.e., a solution of \eqref{eq:equilibria} satisfying $I_0=0$. The fourth equation in \eqref{eq:equilibria} gives $E_0=0$. The fifth, sixth, seventh, first, and second equations then give, respectively, $Q_0$, $H_0$, $R_0$, $S_0$, and $V_0$. In explicit form,
$$S_0=k_6,\,\,\,\,V_0=\frac{u_1k_6}{\mu},\,\,\,\,E_0=0,\,\,\,I_0=0,\,\,\,\,Q_0=\frac{\lambda'}{k_3},\,\,\,\,H_0=\frac{\tau\lambda'}{k_3 k_4},\,\,\,\,R_0=\frac{\lambda'\left(\kappa k_4+\varphi\tau\right)}{k_3 k_4 k_5},$$
where
\begin{eqnarray}\label{eq:kis}
&\displaystyle k_1=\theta+u_3+\mu,\,\,\,k_2=\gamma+u_4+u_5+\mu+\mu',\,\,\,k_3=\kappa+\tau+\mu,\,\,\,k_4=\varphi+\mu+\mu',&\nonumber\\[0.05cm]
&\displaystyle k_5=\mu+\alpha,\,\,\,\text{and}\,\,\,k_6=\frac{\lambda}{u_1+\mu}+\frac{\alpha\lambda'\kappa}{(u_1+\mu)k_3 k_5}+\frac{\alpha\lambda'\varphi\tau}{(u_1+\mu)k_3 k_4 k_5}.&
\end{eqnarray}
The disease-free equilibrium $\mathbf{e}_0$ thus exists ---since all its components are positive--- and is unique, for every set of parameter values.

\sloppy Let us now compute the model's basic reproduction number, using the so-called \textit{next-generation matrix method} \cite[page 33]{DriesscheWatmough}, taking into account as infected compartments those of exposed, infected, quarantined, and hospitalised individuals, whose numbers evolve at the rates given by the third, fourth, fifth, and sixth equations of the model \eqref{eq:model}. Letting $\left(X_1,X_2,X_3,X_4\right):=\left(E,I,Q,H\right)$, we first define
$$\mathcal{F}_1 := \left(1-u_2\right)\beta SX_2+(1-\delta)\beta VX_2,\quad\mathcal{F}_2:=0,\quad \mathcal{F}_3:=0,\quad \mathcal{F}_4:=0,$$
and
\begin{align*}
\mathcal{V}_1 &:= \theta X_1+ u_1X_1 + \mu X_1,\\
\mathcal{V}_2 &:= -\theta X_1 +\gamma X_2 + u_4 X_2 + u_5 X_2 +\mu X_2 +\mu'X_2,\\
\mathcal{V}_3 &:= -\lambda'-u_3X_1-u_4X_2+\kappa X_3+\tau X_3+\mu X_3,\\
\mathcal{V}_4 &:= -\tau X_3-u_5 X_2 +\varphi X_4+\mu X_4+\mu' X_4.
\end{align*}
Next, we define the $4\times 4$ matrices
$$\mathbf{F}:=\left(\begin{array}{c}
\displaystyle\frac{\partial\mathcal{F}_i}{\partial X_j}\left(\mathbf{e}_0\right)
\end{array}\right)=\left(\begin{array}{cccc}0&\left(1-u_2\right)\beta k_6+\left(1-\delta\right)\beta u_1k_6/\mu& 0&0\\0&0&0&0\\0&0&0&0\\0&0&0&0\end{array}\right).$$
and
$$\mathbf{V}:=\left(\begin{array}{c}
\displaystyle\frac{\partial\mathcal{V}_i}{\partial X_j}\left(\mathbf{e}_0\right)
\end{array}\right)=\left(\begin{array}{cccc}
k_1 & 0 & 0 & 0\\
-\theta & k_2 & 0 & 0\\
-u_3 & -u_4 & k_3 & 0\\
0 & -u_5 & -\tau &k_4
\end{array}\right).$$
The basic reproduction number of the model \eqref{eq:model} is the spectral radius of the model's next-generation matrix $\mathbf{F}\mathbf{V}^{-1}$:
\begin{equation}\label{eq:R0}
\cR_0:=\rho\left(\mathbf{F}\mathbf{V}^{-1}\right)=\frac{\theta\beta k_6 \bigl(\mu(1-u_2)+(1-\delta)u_1\bigr)}{k_1 k_2 \mu}.
\end{equation}
Therefore, the basic reproduction number $\cR_0$ grows only sublinearly with the vaccination rate $u_1$. This means that, for the eradication of COVID-19, it is not advisable to strive \textit{only} towards a high vaccination rate; indeed, many of the countries with high percentages of citizens vaccinated \cite{WolfMatthewsAlas} retain their pandemic status. Instead, since $\cR_0$ grows linearly with the mobility intervention rate $u_2$, and with the vaccine efficacy $\delta$, these parameters deserve more attention. In subsection \ref{subsec:sensitivity}, we shall confirm quantitatively that this is the case, i.e., that these are the parameters upon which $\cR_0$ depends most sensitively in the cases of $\cR_0<1$ and $\cR_0>1$, respectively.

Direct computation shows that the characteristic polynomial of the Jacobian matrix of the model \eqref{eq:model} evaluated at $\mathbf{e}_0$ is given by
$$P(x)=\left(x+\mu\right)\left(x+u_1+\mu\right)\left(x+k_3\right)\left(x+k_4\right)\left(x+k_5\right)\left(x^2+bx+c\right),$$
where
\begin{align*}
b &:= 2\mu + \gamma + \theta + u_3+u_4+u_5 +\mu',\\
c &:= \mu^2 + \left(\gamma + \theta+u_3+u_4+u_5+\mu'\right)\mu + \left(\left(u_2-1\right)\beta k_6+\gamma+u_4+u_5+\mu'\right)\theta\\
&\phantom{:=} + u_3\left(\gamma+u_4+u_5+\mu'\right)+\frac{\beta\theta k_6u_1\left(\delta-1\right)}{\mu},
\end{align*}
meaning that $-\mu$, $-u_1-\mu$, $-k_3$, $-k_4$, and $-k_5$ are five negative roots of $P(x)$. Therefore, the equilibrium $\mathbf{e}_0$ is locally asymptotically stable if the other two roots $x_1$ and $x_2$, i.e., those of $x^2+bx+c$, have negative real parts \cite[Theorem 4.6(a)]{Robinson}, and is unstable if at least one of $x_1$ and $x_2$ have a positive real part. Since $b>0$, the former holds if $c>0$ (by the Routh-Hurwitz criterion \cite[section 4.5]{Allen}), while the latter holds if $c<0$ (in which case $x_1$ and $x_2$ are real and have opposite signs). Direct computation shows that $c>0$ is equivalent to $\cR_0<1$, while $c<0$ is equivalent to $\cR_0>1$. This proves the following theorem.\smallskip

\begin{theorem}\label{thm:DFE}
For every set of parameter values, the model \eqref{eq:model} has a unique disease-free equilibrium, which is locally asymptotically stable if $\cR_0 < 1$, and unstable if $\cR_0>1$.
\end{theorem}\smallskip

\subsection{Endemic equilibria}\label{subsec:EE}

Let us now seek all equilibria $\textbf{e}_n=\left(S_n,V_n,E_n,I_n,Q_n,H_n,R_n\right)$ with $I_n\neq 0$, $n\in\mathbb{N}$. The fifth, sixth, seventh, first, and second equations in \eqref{eq:equilibria} give, respectively,
\begin{align}
Q_n &= \frac{\lambda' +u_3E_n+u_4I_n}{k_3},\label{eq:Qn}\\
H_n &= \frac{1}{k_4}\left(\frac{\tau}{k_3}\left(\lambda'+u_3E_n+u_4I_n\right)+u_5I_n\right),\label{eq:Hn}\\
R_n &= \frac{1}{k_5}\left(\gamma I_n + \frac{\kappa}{k_3}\left(\lambda'+u_3E_n+u_4I_n\right)+\frac{\varphi}{k_4}\left(\frac{\tau}{k_3}\left(\lambda'+u_3E_n+u_4I_n\right)+u_5I_n\right)\right),\label{eq:Rn}\\
S_n &= \frac{1}{\mu+u_1+\left(1-u_2\right)\beta I_n}\left(\lambda+\frac{\alpha}{k_5}\left(\gamma I_n + \frac{\kappa}{k_3}\left(\lambda'+u_3E_n+u_4I_n\right)\right.\right.\nonumber\\
&\phantom{=} \left.\left.+\,\frac{\varphi}{k_4}\left(\frac{\tau}{k_3}\left(\lambda'+u_3E_n+u_4I_n\right)+u_5I_n\right)\right)\right)\label{eq:Sn}\\
V_n &= \frac{u_1 S_n}{\left(1-\delta\right)\beta I_n+\mu},\label{eq:Vn}
\end{align}
while the third and fourth equation give
\begin{equation}\label{eq:En}
E_n=\frac{1}{k_1}\left(\left(1-u_2\right)\beta S_n I_n+\left(1-\delta\right)\beta V_n I_n\right)\qquad\text{and}\qquad E_n=\frac{k_2 I_n}{\theta}.
\end{equation}
Equating the two equations in \eqref{eq:En} yields
$$V_n=\frac{\beta\theta u_2 S_n-\beta \theta S_n +k_1k_2}{\theta\left(1-\delta\right)\beta}.$$
Equating this and \eqref{eq:Vn} yields the following expression of $S_n$ as a function of $I_n$:
\begin{equation}\label{eq:Sn1}
S_n=\frac{k_1k_2\left(\mu + \beta I_n - \beta\delta I_n\right)}{\beta\theta\left(\beta\delta u_2 I_n-\beta\delta I_n - \beta u_2 I_n + \beta I_n - u_1\delta - \mu u_2 + \mu + u_1\right)}.
\end{equation}
On the other hand, substituting the second equation in \eqref{eq:En} and into \eqref{eq:Sn} yields another expression of $S_n$ as a function of $I_n$:
\begin{align}\label{eq:Sn2}
S_n&= \frac{1}{\mu+u_1+\left(1-u_2\right)\beta I_n}\left(\lambda+\frac{\alpha}{k_5}\left(\gamma I_n + \frac{\kappa}{k_3}\left(\lambda'+\frac{k_2 u_3I_n}{\theta}+u_4I_n\right)\right.\right.\nonumber\\
&\phantom{=} \left.\left.+\,\frac{\varphi}{k_4}\left(\frac{\tau}{k_3}\left(\lambda'+\frac{k_2u_3 I_n}{\theta}+u_4I_n\right)+u_5I_n\right)\right)\right).
\end{align}
Equating \eqref{eq:Sn1} and \eqref{eq:Sn2}, one finds that the values of $I_n$ are the roots of the quadratic polynomial
\begin{equation}\label{eq:polynomial}
d{I_n}^2 + e I_n + f,
\end{equation}
where
\begin{align*}
d &:= \left(1-\delta\right)\theta\left(1-u_2\right)\beta^2\left(\left(\left(\left(\gamma k_3+\kappa u_4\right)k_4 + \varphi\left(k_3u_5+\tau u_4\right)\right)\theta+ u_3k_2\left(k_4\kappa+\tau\varphi\right)\right)\alpha\right.\\
&\phantom{:=}\left.-\,k_1k_2k_3k_4k_5\right),\\
e &:= \theta\left(\left(\left(\left(\gamma\left(\left(1-u_2\right)\mu+u_1\left(1-\delta\right)\right)k_3-\left(u_4\left(u_2-1\right)\mu+\left(u_1u_4-\beta\lambda'\left(u_2-1\right)\right)\left(\delta-1\right)\right)\kappa\right)k_4\right.\right.\right.\\
&\phantom{:=}\left.\left.\left.-\, \varphi\left(u_5\left(\left(u_2-1\right)\mu+u_1\left(\delta-1\right)\right)k_3+\tau\left(u_4\left(u_2-1\right)\mu+\left(u_1u_4-\beta\lambda'\left(u_2-1\right)\right)\left(\delta-1\right)\right)\right)\right)\theta\right.\right.\\
&\phantom{:=}\left.\left.-\,u_3k_2\left(k_4\kappa+\tau\varphi\right)\left(\left(u_2-1\right)\mu+u_1\left(\delta-1\right)\right)\right)\alpha+k_3k_4k_5\left(\beta\lambda\left(u_2-1\right)\left(\delta-1\right)\theta\right.\right.\\
&\phantom{:=}\left.\left.\,+\left(\left(\delta-u_2-2\right)\mu+u_1\left(\delta-1\right)\right)k_1k_2\right)\right)\beta,\\
f &:= \theta\left(\left(\left(1-u_2\right)\mu+u_1\left(1-\delta\right)\right)\left(\left(\alpha\kappa\lambda'+k_3k_5\lambda\right)k_4+\alpha\varphi\lambda'\tau\right)\beta\theta - k_1k_2k_3k_4k_5\mu\left(\mu+u_1\right)\right).
\end{align*}
Direct computation shows that the condition $\cR_0>1$ is equivalent to $f/d<0$. If this holds, then the values of $I_n$ are real (since $fd<0$) and have opposite signs: $I_1>0$ and $I_2<0$, say. Substituting $I_1$ into \eqref{eq:Qn}, \eqref{eq:Hn}, \eqref{eq:Rn}, \eqref{eq:Sn}, \eqref{eq:Vn}, and \eqref{eq:En}, one obtains a unique endemic equilibrium $\mathbf{e}_1=\left(S_1,V_1,E_1,I_1,Q_1,H_1,R_1\right)$ of the model \eqref{eq:model}, with all components positive. We have therefore proved the following theorem.\smallskip

\begin{theorem}\label{thm:EE}
If $\cR_0>1$, then the model \eqref{eq:model} has a unique positive endemic equilibrium.
\end{theorem}\smallskip

In the case of $\cR_0<1$, we have $f/d>0$, and so no immediate conclusion can be drawn on whether the values of $I_n$ are real. Obtaining an analytic expression for the polynomial's discriminant $e^2-4df$ requires tedious computations, let alone examining its non-negativity. The same situation is faced as we attempt to characterise the stability of these endemic equilibria in the case of their existence, since the analytic expression of $I_n$ given by the quadratic formula is already complicated. This forces us to migrate from analytical to numerical techniques. Preliminary numerical experiments show that, for the parameter values shown in Table \ref{tab:parameters}, $\delta=0.653$, $u_1=10^{-8}$, and $u_2=0.93$, in which case $\cR_0$ is less than but very close to $1$ (cf.\ first case in subsection \ref{subsec:sensitivity}), the polynomial \eqref{eq:polynomial} has two negative real roots, suggesting that the bifurcation occurring at $\cR_0=1$ is a forward transcritical bifurcation \cite[subsection 3.4.3]{Martcheva}. Subsequently, replacing $u_1$ and $u_2$ with $0.4$ and $0.278$, respectively, we observe that $\cR_0$ is greater than $1$ (cf.\ second case in subsection \ref{subsec:sensitivity}) and that a solution of the model \eqref{eq:model} converges towards the unique positive endemic equilibrium guaranteed to exist by Theorem \ref{thm:EE} (cf.\ red graph in Figure \ref{fig:delta}), suggesting that this equilibrium is stable in the case of its existence.\smallskip

\section{Formulation of intervention strategies}\label{section:analysis}

We have mentioned the five concrete forms of intervention presently realised by the Indonesian government to strive towards a new normal: vaccinations, social restrictions, tracings, testings, and treatments. Now, we are ready to exploit the model \eqref{eq:model} to formulate strategies for realising these forms of intervention, in order to optimise their impact. This will be achieved via a two-stage analysis ---which is both numerical and interpretative--- of the model's basic reproduction number $\cR_0$. The first stage is the preliminary analysis, where we identify the set of parameter values corresponding to the \textit{disease-free region} ---that in which $\cR_0<1$--- and its realisability, in various epidemic scenarios. The results, as we shall see, point towards the necessity of vaccinations, and, more importantly, the vaccine efficacy, as important keys to achieve a new normal. The second stage consists in an analysis of the sensitivity of $\cR_0$ with respect to each parameter. For each of the two cases $\cR_0<1$ and $\cR_0>1$, we choose a set of parameter values and compute the sensitivity indices of $\cR_0$ with respect to each parameter, using the results to rank the above five intervention forms in order of significance.

\subsection{Preliminary analysis}\label{subsec:R0}

The five intervention forms are not all equal in the current degree of realisation: tracings, testings, and treatments ---the so-called ``3Ts''--- are reportedly suboptimal \cite{MuthiarinyMurti}, while vaccinations and social restrictions seem to be given primary attention \cite{WHOI19Jan2022,COVID19VaccineTracker,Afifa,Hutasoit}. Accordingly, in this stage of our analysis, let us  assume that the parameters $u_3$, $u_4$, and $u_5$, which represent the rates of the 3Ts, have fixed values. Furthermore, let us fix the values of all parameters except $\delta$, $u_1$, and $u_2$.

As noted in section \ref{section:Introduction}, $u_2=0$ represents normal mobility, while $u_2=1$ represents a total lockdown. The values of $u_2$ which represent social restrictions of level 1, 2, 3, and 4 can be estimated in the following way. First, we deal with level 1 social restrictions, which, as detailed in \cite{SaptoyoKurniawan}, consist of the following regulations:
\begin{enumerate}
\item[(i)] businesses in non-essential sectors are to implement the work-from-office policy at up to $p_1=75\%$ capacity;
\item[(ii)] businesses in essential sectors are to implement the work-from-office policy at up to $p_2=100\%$ capacity;
\item[(iii)] daily-need shops are to operate at up to $p_3=75\%$ capacity;
\item[(iv)] non-daily-need shops are to operate at up to $p_4=75\%$ capacity;
\item[(v)] malls and shopping centres are to operate at up to $p_5=75\%$ capacity;
\item[(vi)] roadside stalls and street vendors are to operate at up to $p_6=75\%$ capacity;
\item[(vii)] restaurants are to operate at up to $p_7=75\%$ capacity;
\item[(viii)] educational activities are to be carried out $p_8=50\%$ onsite and 50\% online;
\item[(ix)] places of worship are to operate at up to $p_9=50\%$ capacity.
\end{enumerate}
We estimate $u_2$ as the average percentage of restrictions in the case of level 1 social restrictions according to the above data: $u_2=1-\left(\sum_{i=1}^9 p_i\right)/9=0.278$. In a similar way, we obtain the following values of $u_2$ representing social restrictions of level 2, 3, and 4: $0.389$, $0.694$, and $0.861$, respectively. The values of $\delta$, on the other hand, will be chosen in view of the efficacies of the actual COVID-19 vaccines \cite{JakartaGlobe,JaraEtAl,MascellinoTimoteoAngelisOliva}.

It is apparent from \eqref{eq:R0} and the definition of the $k_i$s in \eqref{eq:kis} that, for any given $\delta$, the graph $\cR_0=1$ on the $u_1u_2$-plane is a straight line, the abscissa and ordinate intercepts being, respectively,
$$\ell_1:=\frac{\mu\left(\alpha\beta k_4\kappa\theta\lambda'+\alpha\beta\tau\theta\varphi\lambda'+\beta k_3k_4k_5\lambda\theta-k_1k_2k_3k_4k_5\mu\right)}{\alpha\beta\delta k_4 \kappa\theta\lambda'+\alpha\beta\delta\tau\theta\varphi\lambda'+\beta\delta k_3k_4k_5\lambda\theta-\alpha\beta k_4 \kappa\theta\lambda'-\alpha\beta\tau\theta\varphi\lambda'-\beta k_3k_4k_5\lambda\theta+k_1k_2k_3k_4k_5\mu}$$
and
$$\ell_2:=\frac{\left(\kappa\theta\beta\alpha\lambda'+k_3 k_5\left(\beta\lambda\theta-k_1k_2\mu\right)\right)k_4+\alpha\beta\tau\theta\varphi\lambda'}{\left(\left(\alpha\kappa\lambda'+k_3k_5\lambda\right)k_4+\varphi\lambda'\tau\alpha\right)\beta\theta}.$$
Notice that $\ell_2$ is independent of $\delta$, since so are the $k_i$s, by \eqref{eq:kis}. Moreover, a direct computation shows that the denominator of $\ell_1$ is equal to zero if and only if $\delta=\ell_2$. Furthermore, we let
$$\ell_3:=\frac{k_1k_2k_3k_4k_5\mu^2}{\alpha\beta k_4\kappa\theta\lambda'+\alpha\beta\tau\theta\varphi\lambda
+\beta k_3k_4k_5\lambda\theta-k_1k_2k_3k_4k_5\mu-\alpha\beta\delta k_4\kappa\theta\lambda'-\alpha\beta\delta\tau\theta\varphi\lambda'-\beta\delta k_3k_4k_5\lambda\theta}$$
be the abscissa of the point of ordinate $1$ on the line. Substituting the values shown in Table \ref{tab:parameters} of all parameters except $\delta$, $u_1$, and $u_2$, one obtains
\begin{align*}
\ell_1&=\frac{0.0000139049}{0.3549600264\delta-0.3298930314},\\
\ell_2&= 0.9293807942,\\
\ell_3&=\frac{0.0000042150}{1.3160453862-1.4160453862\delta}.
\end{align*}

\begin{figure}\centering
\begin{tabular}{ccc}
(1)&(2)&(3)\\
\scalebox{0.85}{\begin{tikzpicture}
\pgfset{declare function={f(\x)=.9293807946+6557.134355*\x;}}

\begin{axis}[
	xmin=-0.000155,
	xmax=0.000032,
	ymin=-0.11,
	ymax=1.1,
	xtick={-0.0001417358,0.1076982742e-4},
	xticklabels={{\footnotesize $\ell_1$},{\footnotesize $\ell_3$}},
	ytick={0.9293807953,1},
	yticklabels={{\footnotesize $\ell_2$},{\footnotesize $1$}},
	axis lines=middle,
	axis on top=true,
	samples=100,
	xlabel=$u_1$,
	ylabel=$u_2$,
	x axis line style=<->,
	y axis line style=<->,
	width=6.75cm,
	height=6.75cm,
	scaled x ticks=false
]
\draw[dashed] (axis cs:0,1) -- (axis cs:0.00003,1);
\draw[dashed] (axis cs:0.1076982742e-4,0) -- (axis cs:0.1076982742e-4,1);
\addplot [very thick,domain=-0.0001417358:0.1076982742e-4] {f(x)};
\node[above left,xshift=5pt] at (axis cs:-0.00007086790000,.4646904528) {$\cR_0=1$};
\fill[pattern=north west lines] (axis cs:0,1) -- (axis cs:0,0.9293807953) -- (axis cs:0.1076982742e-4,1) -- cycle;
\end{axis}
\end{tikzpicture}}&\scalebox{0.85}{\begin{tikzpicture}
\pgfset{declare function={f(\x)=(1.96830308683889+1.16576998454192*10^5*\x)/(23725.*\x+1.);}}
\pgfset{declare function={g(\x)=(0.708022693107513e-4+1.16576998454192*10^5*\x)/(23725.*\x+1.);}}

\begin{axis}[
	xmin=0,
	xmax=0.0004,
	ymin=0,
	ymax=5,
	xtick={0,0.0002,0.0004},
	xticklabels={0,0.0002,0.0004},
	ytick={1,2,3,4,5},
	yticklabels={1,2,3,4,5},
	axis on top=true,
	samples=100,
	xlabel=$u_1$,
	ylabel=$\cR_0$,
	width=5.5cm,
	height=5.5cm,
	ylabel near ticks,
	ticklabel style = {font=\footnotesize}
]
\addplot [red,very thick,domain=0:0.0004] {f(x)};
\addplot [blue,very thick,domain=0:0.0004] {g(x)};
\end{axis}
\end{tikzpicture}}&\scalebox{0.85}{\begin{tikzpicture}
\pgfset{declare function={f(\x)=4.914651762-0.1491987553e-2*\x;}}
\pgfset{declare function={g(\x)=13.179876057142695166-12.658803318493943770*\x;}}

\begin{axis}[
	xmin=0,
	xmax=1,
	ymin=0,
	ymax=13.5,
	xtick={0,0.25,0.50,0.75,1},
	xticklabels={0,0.25,0.50,0.75,1},
	ytick={1,5,9,13},
	yticklabels={1,5,9,13},
	axis on top=true,
	samples=100,
	xlabel=$u_2$,
	ylabel=$\cR_0$,
	width=5.5cm,
	height=5.5cm,
	ylabel near ticks,
	ticklabel style = {font=\footnotesize}
]
\addplot [red,very thick,domain=0:1] {f(x)};
\addplot [blue,very thick,domain=0:1] {g(x)};
\end{axis}
\end{tikzpicture}}
\end{tabular}
\caption{\label{fig:R0experiment1}(1) Plot of the line $\cR_0=1$ on the $u_1u_2$-plane in the case of $\delta=0.653$, with the feasible disease-free region shaded; (2) plot of $\cR_0$ as a function of $u_1$ in the case of $\delta=0.653$, for $u_2=0.861$ (red) and for $u_2=0.999995$ (blue); (3) plot of $\cR_0$ as a function of $u_2$ in the case of $\delta=0.653$, for $u_1=0.4$ (red) and for $u_1=0.000005$ (blue).}
\end{figure}
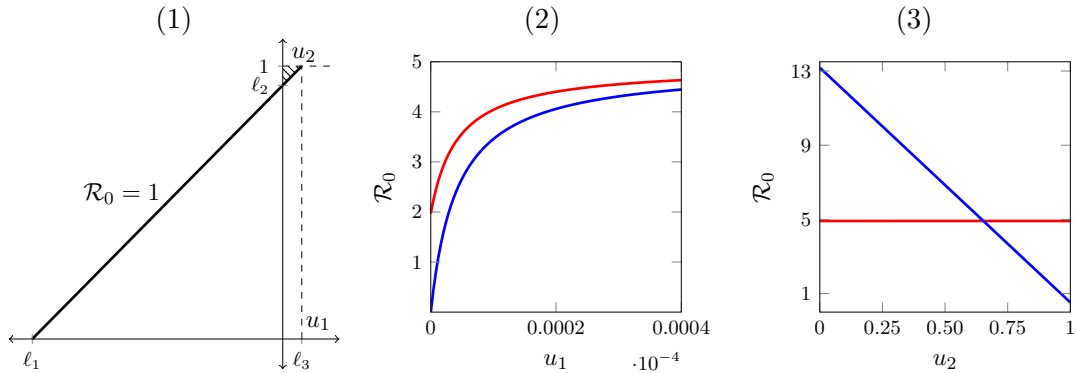

Recently, the vaccine most sought-after in the country is reportedly Sinovac's Coronavac \cite{Bona}, which has demonstrated a $65.3\%$ efficacy \cite{JakartaGlobe,JaraEtAl}. For this value of $\delta$,
we have $\ell_1=-0.0001417358$ and $\ell_3=0.0000107698$, and the line $\cR_0$ on the $u_1u_2$-plane is plotted in Figure \ref{fig:R0experiment1} (1). The shaded region is the feasible disease-free region, i.e., the region $\left\{\left(u_1,u_2\right)\in[0,1]^2:\cR_0<1\right\}$. Therefore, according to our model, using a vaccine with only a $65.3\%$ efficacy, the pandemic can only be resolved if $u_1\leqslant\ell_3=0.0000107698$, i.e., the vaccination rate is made extremely low, and $u_2\geqslant\ell_2=0.9293807942$, i.e., a near-lockdown policy is implemented. The latter is uncompromisable: even level 4 social restrictions are insufficient; see the red curve in Figure \ref{fig:R0experiment1} (2). Likewise, if the vaccination rate is increased even only to a moderate level, say $u_1=0.4$ (Table \ref{tab:parameters}), then the policy of raising the level of social restrictions becomes insignificant: such a policy suppresses the endemic-valued basic reproduction number ---and thus the number of daily new cases--- only insignificantly; see Figure \ref{fig:R0experiment1} (3). A major reason for this is that, for vaccinated citizens, social restrictions are waived\footnote{Notice in the model's compartment diagram (Figure \ref{fig:presentmodel}) that the mobility restriction factor $1-u_2$ is present in the transition rate from compartment S to compartment E, but not in the transition rate from compartment V to compartment E.}, allowing them to travel, visit public places, etc.\ more unrestrictedly than unvaccinated citizens, bringing about a high risk in the case of high vaccination rate but low vaccine efficacy. We thus find it unsurprising that the omicron wave remained unavoidable despite the notable progress of the country's vaccination programme (see section \ref{section:Introduction}).

\begin{figure}\centering
\begin{tabular}{ccc}
(1)&(2)&(3)\\
\scalebox{0.85}{\begin{tikzpicture}
\pgfset{declare function={f(\x)=.9293807946+697.0593528*\x;}}

\begin{axis}[
	xmin=-0.0015,
	xmax=0.00031,
	ymin=-0.11,
	ymax=1.1,
	xtick={-0.0013332879,0.0001013102},
	xticklabels={{\footnotesize $\ell_1$},{\footnotesize $\ell_3$}},
	ytick={0.9293807953,1},
	yticklabels={{\footnotesize $\ell_2$},{\footnotesize $1$}},
	axis lines=middle,
	axis on top=true,
	samples=100,
	xlabel=$u_1$,
	ylabel=$u_2$,
	x axis line style=<->,
	y axis line style=<->,
	width=6.75cm,
	height=6.75cm,
	scaled x ticks=false
]
\draw[dashed] (axis cs:0,1) -- (axis cs:0.00029,1);
\draw[dashed] (axis cs:0.1013101755e-3,0) -- (axis cs:0.1013101755e-3,1);
\addplot [very thick,domain=-0.0013332879:0.0001013102] {f(x)};
\node[above left,xshift=5pt] at (axis cs:-0.6666439500e-3,.4646903943) {$\cR_0=1$};
\fill[pattern=north west lines] (axis cs:0,1) -- (axis cs:0,0.9293807953) -- (axis cs:0.1013101755e-3,1) -- cycle;
\end{axis}
\end{tikzpicture}}&\scalebox{0.85}{\begin{tikzpicture}
\pgfset{declare function={f(\x)=(49.2075771709722+8.3989191969879*10^5*\x)/(5.93125*10^5*\x+25.);}}
\pgfset{declare function={g(\x)=(0.708022693107512e-4+33595.6767879516*\x)/(23725.*\x+1.);}}

\begin{axis}[
	xmin=0,
	xmax=0.0004,
	ymin=0,
	ymax=2,
	xtick={0,0.0002,0.0004},
	xticklabels={0,0.0002,0.0004},
	ytick={0.5,1.0,1.5,2},
	yticklabels={0.5,1.0,1.5,2},
	axis on top=true,
	samples=100,
	xlabel=$u_1$,
	ylabel=$\cR_0$,
	width=5.5cm,
	height=5.5cm,
	ylabel near ticks,
	ticklabel style = {font=\footnotesize}
]
\addplot [red,very thick,domain=0:0.0004] {f(x)};
\addplot [blue,very thick,domain=0:0.0004] {g(x)};
\end{axis}
\end{tikzpicture}}&\scalebox{0.85}{\begin{tikzpicture}
\pgfset{declare function={f(\x)=1.417388176 - 0.001491987553*\x;}} 
\pgfset{declare function={g(\x)=1.416642217 - 0.0006631443963*\x;}} 

\begin{axis}[
	xmin=0,
	xmax=1,
	ymin=1.415896188,
	ymax=1.417388176,
	xtick={0,0.25,0.50,0.75,1},
	xticklabels={0,0.25,0.50,0.75,1},
	ytick={1.4160,1.4164,1.4168,1.4172},
	yticklabels={1.4160,1.4164,1.4168,1.4172},
	axis on top=true,
	samples=100,
	xlabel=$u_2$,
	ylabel=$\cR_0$,
	width=5.5cm,
	height=5.5cm,
	ylabel near ticks,
	ticklabel style = {font=\footnotesize}
]
\addplot [red,very thick,domain=0:1] {f(x)};
\addplot [blue,very thick,domain=0:1] {g(x)};
\end{axis}
\end{tikzpicture}}
\end{tabular}
\caption{\label{fig:R0experiment2}(1) Plot of the line $\cR_0=1$ on the $u_1u_2$-plane in the case of $\delta=0.9$, with the feasible disease-free region shaded; (2) plot of $\cR_0$ as a function of $u_1$ in the case of $\delta=0.9$, for $u_2=0.861$ (red) and for $u_2=0.999995$ (blue); (3) plot of $\cR_0$ as a function of $u_2$ in the case of $\delta=0.9$, for $u_1=0.4$ (red) and for $u_1=0.9$ (blue).}
\end{figure}
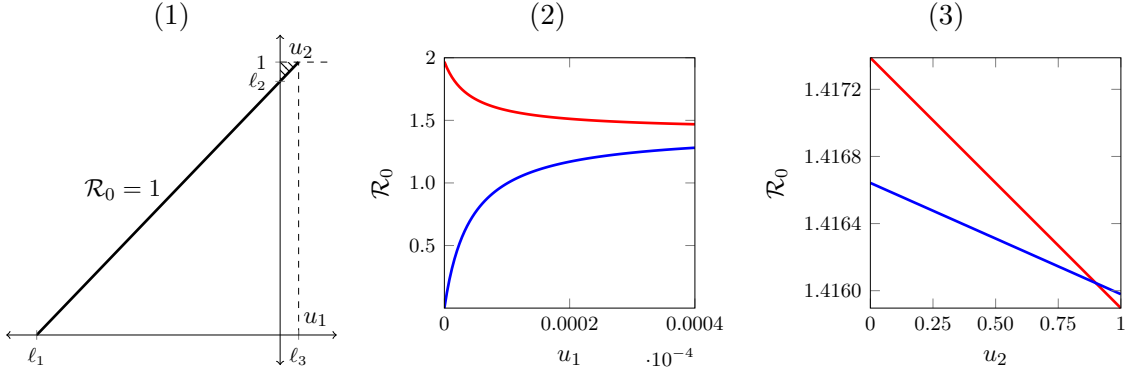

Now let us suppose that the country utilises a vaccine with a higher efficacy: say, $\delta=0.9$. In this case, we have $\ell_1=-0.0013332879$ and $\ell_3=0.0001013102$, and the line $\cR_0$ on the $u_1u_2$-plane is plotted in Figure \ref{fig:R0experiment2} (1). Since no qualitative change is observed here, the message remains the same: the disease's transmission can only be halted with an extremely low vaccination rate and an extremely high level of mobility restrictions. However, comparing Figure \ref{fig:R0experiment1} (2) and Figure \ref{fig:R0experiment2} (2), we see a qualitative change: the red curve, which, in both figures, correspond to level 4 social restrictions, i.e., $u_1=0.861$, has changed monotonicity. Furthermore, comparing the quantitative properties of the red lines in Figure \ref{fig:R0experiment1} (3) and Figure \ref{fig:R0experiment2} (3), both corresponding to the moderate vaccination rate $u_1=0.4$, we can see that the improvement of the vaccine efficacy, from $0.653$ to $0.9$, drastically decreases the value of the basic reproduction number, from above 4 to below 2. We infer therefore that the improvement of the quality of COVID-19 vaccines should take precedence over that of the rate at which vaccinations are carried out.

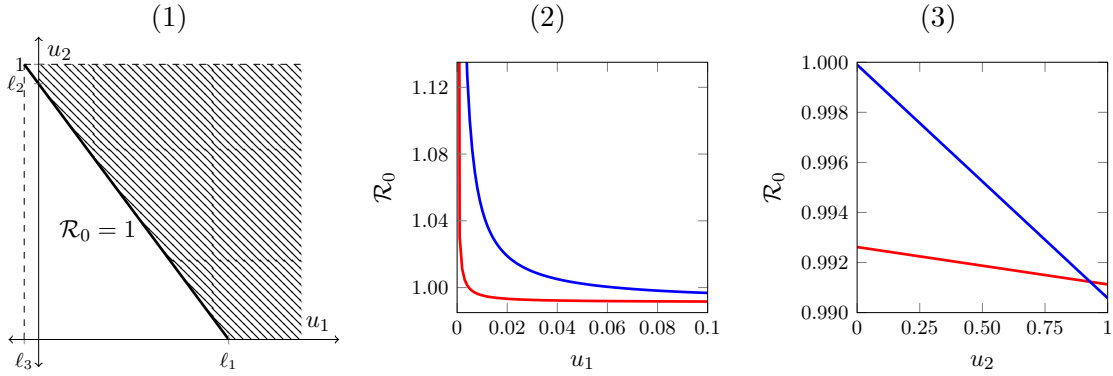
\begin{figure}\centering
\begin{tabular}{ccc}
(1)&(2)&(3)\\
\scalebox{0.85}{\begin{tikzpicture}
\pgfset{declare function={f(\x)=.9293807946-14.69064745*\x;}}

\begin{axis}[
	xmin=-0.01,
	xmax=0.1,
	ymin=-0.1,
	ymax=1.1,
	xtick={-0.0048070860,0.0632634203},
	xticklabels={{\footnotesize $\ell_3$},{\footnotesize $\ell_1$}},
	ytick={0.9293807953,1},
	yticklabels={{\footnotesize $\ell_2$},{\footnotesize $1$}},
	axis lines=middle,
	axis on top=true,
	samples=100,
	xlabel=$u_1$,
	ylabel=$u_2$,
	x axis line style=<->,
	y axis line style=<->,
	width=6.75cm,
	height=6.75cm,
	scaled x ticks=false
]
\draw[dashed] (axis cs:-0.0048070860,1) -- (axis cs:0.0875,1);
\draw[dashed] (axis cs:-0.0048070860,0) -- (axis cs:-0.0048070860,1);
\addplot [very thick,domain=-0.004807085980:0.06326343327] {f(x)};
\node[below left,xshift=5pt] at (axis cs:0.3163171015e-1,0.4646904925) {$\cR_0=1$};
\fill[pattern=north west lines] (axis cs:0,1) -- (axis cs:0,0.9293807953) -- (axis cs:0.0632634203,0) -- (axis cs:0.0875,0) -- (axis cs:0.0875,1) -- cycle;
\end{axis}
\end{tikzpicture}} & \scalebox{0.85}{\begin{tikzpicture}
\pgfset{declare function={f(\x)=(1.96830308683889+23516.9737515661*\x)/(23725.*\x+1.);}}
\pgfset{declare function={g(\x)=(14.1604538621503+23516.9737515661*\x)/(23725.*\x+1.);}}

\begin{axis}[
	xmin=0,
	xmax=0.1,
	ymin=0.985,
	ymax=1.135,
	xtick={0,0.02,0.04,0.06,0.08,0.1},
	xticklabels={0,0.02,0.04,0.06,0.08,0.1},
	ytick={1.00,1.04,1.08,1.12},
	yticklabels={1.00,1.04,1.08,1.12},
	axis on top=true,
	samples=100,
	xlabel=$u_1$,
	ylabel=$\cR_0$,
	width=5.5cm,
	height=5.5cm,
	ylabel near ticks,
	ticklabel style = {font=\footnotesize}
]
\addplot [red,very thick,domain=0:0.1] {f(x)};
\addplot [blue,very thick,domain=0:0.1] {g(x)};
\end{axis}
\end{tikzpicture}} & \scalebox{0.85}{\begin{tikzpicture}
\pgfset{declare function={f(\x)=.99261931877447966557-0.14919875526446410125e-2*\x;}}
\pgfset{declare function={g(\x)=0.99989915358850867629 - 0.0093197669225683084436*\x;}}

\begin{axis}[
	xmin=0,
	xmax=1,
	ymin=0.990,
	ymax=1,
	xtick={0,0.25,0.50,0.75,1},
	xticklabels={0,0.25,0.50,0.75,1},
	ytick={0.990,0.992,0.994,0.996,0.998,1.000},
	yticklabels={0.990,0.992,0.994,0.996,0.998,1.000},
	axis on top=true,
	samples=100,
	xlabel=$u_2$,
	ylabel=$\cR_0$,
	width=5.5cm,
	height=5.5cm,
	ylabel near ticks,
	ticklabel style = {font=\footnotesize}
]
\addplot [red,very thick,domain=0:1] {f(x)};
\addplot [blue,very thick,domain=0:1] {g(x)};
\end{axis}
\end{tikzpicture}}
\end{tabular}
\caption{\label{fig:R0experiment3}(1) Plot of the line $\cR_0=1$ on the $u_1u_2$-plane in the case of $\delta=0.93$, with the feasible disease-free region shaded; (2) plot of $\cR_0$ as a function of $u_1$ in the case of $\delta=0.93$, for $u_2=0.861$ (red) and for $u_2=0$ (blue); (3) plot of $\cR_0$ as a function of $u_2$ in the case of $\delta=0.93$, for $u_1=0.4$ (red) and for $u_1=0.064$ (blue).}
\end{figure}

Let us further increase the vaccine efficacy: $\delta=0.93$. In this case, we have $\ell_1=0.0632634203$ and $\ell_3=-0.0048070860$, and the line $\cR_0$ on the $u_1u_2$-plane is plotted in Figure \ref{fig:R0experiment3} (1). Now, we see a radical qualitative change ---a much desirable one--- from Figure \ref{fig:R0experiment3} (1): the line's slope is now negative, and a disease-free state can be achieved even with a complete removal of social restrictions, i.e., $u_2=0$, and a very low vaccination rate, i.e., any low value of $u_1$ satisfying $u_1\geqslant \ell_1=0.0632634203$, say, $u_1=0.064$; see Figure \ref{fig:R0experiment3} (2) and (3).

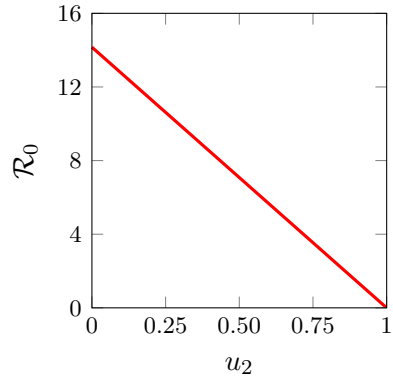
\begin{figure}\centering
\begin{tikzpicture}
\pgfset{declare function={f(\x)=14.16045386-14.16045386*\x;}}

\begin{axis}[
	xmin=0,
	xmax=1,
	ymin=0,
	ymax=16,
	xtick={0,0.25,0.50,0.75,1},
	xticklabels={0,0.25,0.50,0.75,1},
	ytick={0,4,8,12,16},
	yticklabels={0,4,8,12,16},
	axis on top=true,
	samples=100,
	xlabel=$u_2$,
	ylabel=$\cR_0$,
	width=5.5cm,
	height=5.5cm,
	ylabel near ticks,
	ticklabel style = {font=\footnotesize}
]
\addplot [red,very thick,domain=0:1] {f(x)};
\end{axis}
\end{tikzpicture}
\caption{\label{fig:R0experiment4}Plot of $\cR_0$ as a function of $u_2$ in the case of $u_1=0$.}
\end{figure}

A follow-up question naturally arises: is it possible to achieve a disease-free state without vaccinations, i.e., with $u_1=0$? In this case, the expression \eqref{eq:R0} for $\cR_0$ is independent of $\delta$:
$$\cR_0=\frac{\left(1-u_2\right)\beta\left(\alpha k_4 \kappa \lambda'+\alpha\tau\varphi\lambda'+k_3k_4k_5\lambda\right)\theta}{\mu k_1 k_2 k_3 k_4 k_5}=14.1604538645-14.1604538645u_2;$$
a plot of $\cR_0$ versus $u_2$ is shown in Figure \ref{fig:R0experiment4}: a line with a negative slope which is rather large in absolute value. Therefore, in absence of vaccinations, raising the level of social restrictions suppresses the basic reproduction number significantly. This justifies the effectiveness of the government's social restriction policies prior to the commencement of the vaccination programme. Notice that, near $1$, the values of $\cR_0$ as shown Figure \ref{fig:R0experiment3} are lower than those shown in Figure \ref{fig:R0experiment1} (right); that is, in the cases where social restrictions are imposed on level 3 or 4, it is better not to administer vaccines than to administer vaccines with low efficacy, which is again unsurprising for the reason mentioned earlier: restrictions are waived for vaccinated citizens. Nevertheless, in absence of vaccinations, even level 4 social restrictions are not sufficient to bring the country to the disease-free state. The latter requires, again, a near-lockdown policy: $u_2\geqslant \ell_2=0.9293807942$.

Let us summarise the recommended strategies arising from this first-stage analysis. Firstly, if a lockdown is undesirable, it is necessary to administer vaccinations. However, one should strive not primarily towards the increase of the rate at which they are administered, but towards the use of high-efficacy vaccines\footnote{The severity of the vaccines' side effects may also need to be considered for the acceptability of this strategy; alternatives include the utilisation of more than one type of vaccines with different levels of efficacy and side effects, the analysis of which requires a modification of the model \eqref{eq:model}; see section \ref{section:conclusions}.}, such as Pfizer-BioNTech or Moderna \cite{MascellinoTimoteoAngelisOliva}. Secondly, it is necessary to set out and implement an appropriate level of social restrictions to vaccinated citizens\footnote{This motivates another modification of the model \eqref{eq:model}; see section \ref{section:conclusions}.}, especially those who received vaccines with limited efficacy and/or have only been vaccinated partially.\smallskip

\subsection{Sensitivity analysis}\label{subsec:sensitivity}

Let us now complement the above analysis with a quantitative assessment of the significance of each parameter. We compute the \textit{sensitivity index} \cite{ChitnisHymanCushing} of the basic reproduction number $\cR_0$ with respect to a parameter $p$, i.e.,
\begin{equation}\label{eq:sensitivity}
\Upsilon^{\cR_0}_p:=\frac{\partial \cR_0}{\partial p}\cdot\frac{p}{\cR_0},
\end{equation}
for every $p\in\mathcal{P}$, where
$$\mathcal{P}:=\left\{\lambda,\lambda',\mu,\mu',\beta,\delta,\alpha,\theta,\gamma,\varphi,\kappa,\tau,u_1,u_2,u_3,u_4,u_5\right\}$$
denotes the set of all parameters in the model \eqref{eq:model}, obtaining, e.g., for $p=\delta$ and $p=u_2$,
$$\Upsilon^{\cR_0}_{\delta}=\frac{\delta u_1}{\delta u_1+\mu u_2-\mu-u_1}\qquad\text{and}\qquad\Upsilon^{\cR_0}_{u_2}=\frac{\mu u_2}{\delta u_1+\mu u_2-\mu-u_1},$$
respectively.

We choose two sets of parameter values: one representing a disease-free case ($\cR_0<1$), and another representing an endemic case ($\cR_0>1$). In both sets, the values of $\lambda$, $\lambda'$, $\mu$, $\mu'$, $\beta$, $\alpha$, $\theta$, $\gamma$, $\varphi$, $\kappa$, $\tau$, $u_3$, $u_4$, $u_5$ are as shown in Table \ref{tab:parameters}, and $\delta=0.653$. In the former case, we set $\left(u_1,u_2\right)=\left(10^{-8},0.93\right)$, so that $\cR_0 = 0.9921621498<1$, while in the latter, we set $\left(u_1,u_2\right)=\left(0.4,0.278\right)$, so that $\cR_0=4.9142369856>1$.

In each case, we substitute the parameter values to all sensitivity indices. For each index $\Upsilon^{\cR_0}_p$, the following aspects are essential.
\begin{enumerate}
\item \ul{The sign $\sgn\left(\Upsilon^{\cR_0}_p\right)$}, which is positive (negative) if and only if $\cR_0$ is monotonically increasing (decreasing) with $p$.
\item \ul{The absolute value $\left|\Upsilon^{\cR_0}_p\right|$}, which measures the relative change of $\cR_0$ with respect to $p$: a $P\%$ change of $p$ results in a $\left|\Upsilon^{\cR_0}_p\right|P\%$ change of $\cR_0$. Thus, the higher the value of $\left|\Upsilon^{\cR_0}_p\right|$, the more significant the parameter $p$.
\end{enumerate}
It was our intention to visualise the results using a bidirectional bar chart of the values of $\Upsilon^{\cR_0}_p$ for every $p$ in each case, but the rather unusual distribution of these values, especially in the disease-free case where an extreme outlier is present, makes such a chart ineffective.

\begin{figure}\centering
\begin{tikzpicture}
\begin{axis}[
	xmin=-25,
	xmax=22,
	ymin=-0.5,
	ymax=17.5,
	xtick={-16,-12,-8,-4,0,4,8,12,16},
	ytick=\empty,
	axis on top=true,
	samples=100,
	xlabel=$\sgn\left(\Upsilon^{\cR_0}_p\right)\left(18-r(p)\right)$,
	ylabel=$p$,
	width=15cm,
	height=7cm,
	ylabel near ticks,
	ticklabel style = {font=\footnotesize},
	clip=false
]
\node[above] at (axis cs:-1.5,17.5) {$\cR_0<1$};
\draw(axis cs:0,-0.5) -- (axis cs:0,17.5);
\draw[pattern=north east lines] (axis cs:0,0) rectangle (axis cs:-8,1); 
\node[right] at (axis cs:0,0.5) {\tiny $u_5$};
\node[left] at (axis cs:-8,0.5) {\tiny $-0.1625549344$};

\draw[pattern=north east lines] (axis cs:0,1) rectangle (axis cs:-13,2); 
\node[right] at (axis cs:0,1.5) {\tiny $u_4$};
\node[left] at (axis cs:-13,1.5) {\tiny $-0.5854319365$};

\draw[pattern=north east lines] (axis cs:0,2) rectangle (axis cs:-11,3); 
\node[right] at (axis cs:0,2.5) {\tiny $u_3$};
\node[left] at (axis cs:-11,2.5) {\tiny $-0.5555295385$};

\draw[pattern=north east lines] (axis cs:0,3) rectangle (axis cs:-17,4); 
\node[right] at (axis cs:0,3.5) {\tiny $u_2$};
\node[left] at (axis cs:-17,3.5) {\tiny $-13.2701075492$};

\draw[pattern=north east lines] (axis cs:0,4) rectangle (axis cs:5,5); 
\node[left] at (axis cs:0,4.5) {\tiny $u_1$};
\node[right] at (axis cs:5,4.5) {\tiny $0.0009375069$};

\draw[pattern=north east lines] (axis cs:0,5) rectangle (axis cs:-1,6); 
\node[right] at (axis cs:0,5.5) {\tiny $\tau$};
\node[left] at (axis cs:-1,5.5) {\tiny $-0.0005770690$};

\draw[pattern=north east lines] (axis cs:0,6) rectangle (axis cs:3,7); 
\node[left] at (axis cs:0,6.5) {\tiny $\kappa$};
\node[right] at (axis cs:3,6.5) {\tiny $0.0006557000$};

\draw[pattern=north east lines] (axis cs:0,7) rectangle (axis cs:2,8); 
\node[left] at (axis cs:0,7.5) {\tiny $\varphi$};
\node[right] at (axis cs:2,7.5) {\tiny $0.0006205788$};

\draw[pattern=north east lines] (axis cs:0,8) rectangle (axis cs:-9,9); 
\node[right] at (axis cs:0,8.5) {\tiny $\gamma$};
\node[left] at (axis cs:-9,8.5) {\tiny $-0.1951439788$};

\draw[pattern=north east lines] (axis cs:0,9) rectangle (axis cs:12,10); 
\node[left] at (axis cs:0,9.5) {\tiny $\theta$};
\node[right] at (axis cs:12,9.5) {\tiny $0.5555763692$};

\draw[pattern=north east lines] (axis cs:0,10) rectangle (axis cs:4,11); 
\node[left] at (axis cs:0,10.5) {\tiny $\alpha$};
\node[right] at (axis cs:4,10.5) {\tiny $0.0007836083$};

\draw[pattern=north east lines] (axis cs:0,11) rectangle (axis cs:-6,12); 
\node[right] at (axis cs:0,11.5) {\tiny $\delta$};
\node[left] at (axis cs:-6,11.5) {\tiny $-0.0022106037$};

\draw[pattern=north east lines] (axis cs:0,12) rectangle (axis cs:15,13); 
\node[left] at (axis cs:0,12.5) {\tiny $\beta$};
\node[right] at (axis cs:15,12.5) {\tiny $1.0000000000$};

\draw[pattern=north east lines] (axis cs:0,13) rectangle (axis cs:-7,14); 
\node[right] at (axis cs:0,13.5) {\tiny $\mu'$};
\node[left] at (axis cs:-7,13.5) {\tiny $-0.0574065790$};

\draw[pattern=north east lines] (axis cs:0,14) rectangle (axis cs:-16,15); 
\node[right] at (axis cs:0,14.5) {\tiny $\mu$};
\node[left] at (axis cs:-16,14.5) {\tiny $-1.0019299760$};

\draw[pattern=north east lines] (axis cs:0,15) rectangle (axis cs:10,16); 
\node[left] at (axis cs:0,15.5) {\tiny $\lambda'$};
\node[right] at (axis cs:10,15.5) {\tiny $0.2052857854$};

\draw[pattern=north east lines] (axis cs:0,16) rectangle (axis cs:14,17); 
\node[left] at (axis cs:0,16.5) {\tiny $\lambda$};
\node[right] at (axis cs:14,16.5) {\tiny $0.7947142143$};
\end{axis}
\end{tikzpicture}\bigskip

\begin{tikzpicture}
\begin{axis}[
	xmin=-25,
	xmax=22,
	ymin=-0.5,
	ymax=17.5,
	xtick={-16,-12,-8,-4,0,4,8,12,16},
	ytick=\empty,
	axis on top=true,
	samples=100,
	xlabel=$\sgn\left(\Upsilon^{\cR_0}_p\right)\left(18-r(p)\right)$,
	ylabel=$p$,
	width=15cm,
	height=7cm,
	ylabel near ticks,
	ticklabel style = {font=\footnotesize},
	clip=false
]
\node[above] at (axis cs:-1.5,17.5) {$\cR_0>1$};
\draw(axis cs:0,-0.5) -- (axis cs:0,17.5);
\draw[pattern=north east lines] (axis cs:0,0) rectangle (axis cs:-8,1); 
\node[right] at (axis cs:0,0.5) {\tiny $u_5$};
\node[left] at (axis cs:-8,0.5) {\tiny $-0.1625549344$};

\draw[pattern=north east lines] (axis cs:0,1) rectangle (axis cs:-13,2); 
\node[right] at (axis cs:0,1.5) {\tiny $u_4$};
\node[left] at (axis cs:-13,1.5) {\tiny $-0.5854319365$};

\draw[pattern=north east lines] (axis cs:0,2) rectangle (axis cs:-11,3); 
\node[right] at (axis cs:0,2.5) {\tiny $u_3$};
\node[left] at (axis cs:-11,2.5) {\tiny $-0.5555295385$};

\draw[pattern=north east lines] (axis cs:0,3) rectangle (axis cs:-1,4); 
\node[right] at (axis cs:0,3.5) {\tiny $u_2$};
\node[left] at (axis cs:-1,3.5) {\tiny $-0.0000844022$};

\draw[pattern=north east lines] (axis cs:0,4) rectangle (axis cs:-2,5); 
\node[right] at (axis cs:0,4.5) {\tiny $u_1$};
\node[left] at (axis cs:-2,4.5) {\tiny $-0.0001138399$};

\draw[pattern=north east lines] (axis cs:0,5) rectangle (axis cs:-3,6); 
\node[right] at (axis cs:0,5.5) {\tiny $\tau$};
\node[left] at (axis cs:-3,5.5) {\tiny $-0.0005770690$};

\draw[pattern=north east lines] (axis cs:0,6) rectangle (axis cs:5,7); 
\node[left] at (axis cs:0,6.5) {\tiny $\kappa$};
\node[right] at (axis cs:5,6.5) {\tiny $0.0006557000$};

\draw[pattern=north east lines] (axis cs:0,7) rectangle (axis cs:4,8); 
\node[left] at (axis cs:0,7.5) {\tiny $\varphi$};
\node[right] at (axis cs:4,7.5) {\tiny $0.0006205788$};

\draw[pattern=north east lines] (axis cs:0,8) rectangle (axis cs:-9,9); 
\node[right] at (axis cs:0,8.5) {\tiny $\gamma$};
\node[left] at (axis cs:-9,8.5) {\tiny $-0.1951439788$};

\draw[pattern=north east lines] (axis cs:0,9) rectangle (axis cs:12,10); 
\node[left] at (axis cs:0,9.5) {\tiny $\theta$};
\node[right] at (axis cs:12,9.5) {\tiny $0.5555763692$};

\draw[pattern=north east lines] (axis cs:0,10) rectangle (axis cs:6,11); 
\node[left] at (axis cs:0,10.5) {\tiny $\alpha$};
\node[right] at (axis cs:6,10.5) {\tiny $0.0007836083$};

\draw[pattern=north east lines] (axis cs:0,11) rectangle (axis cs:-17,12); 
\node[right] at (axis cs:0,11.5) {\tiny $\delta$};
\node[left] at (axis cs:-17,11.5) {\tiny $-1.8814318750$};

\draw[pattern=north east lines] (axis cs:0,12) rectangle (axis cs:15,13); 
\node[left] at (axis cs:0,12.5) {\tiny $\beta$};
\node[right] at (axis cs:15,12.5) {\tiny $1.0000000000$};

\draw[pattern=north east lines] (axis cs:0,13) rectangle (axis cs:-7,14); 
\node[right] at (axis cs:0,13.5) {\tiny $\mu'$};
\node[left] at (axis cs:-7,13.5) {\tiny $-0.0574065790$};

\draw[pattern=north east lines] (axis cs:0,14) rectangle (axis cs:-16,15); 
\node[right] at (axis cs:0,14.5) {\tiny $\mu$};
\node[left] at (axis cs:-16,14.5) {\tiny $-1.0019299760$};

\draw[pattern=north east lines] (axis cs:0,15) rectangle (axis cs:10,16); 
\node[left] at (axis cs:0,15.5) {\tiny $\lambda'$};
\node[right] at (axis cs:10,15.5) {\tiny $0.2052857854$};

\draw[pattern=north east lines] (axis cs:0,16) rectangle (axis cs:14,17); 
\node[left] at (axis cs:0,16.5) {\tiny $\lambda$};
\node[right] at (axis cs:14,16.5) {\tiny $0.7947142143$};
\end{axis}
\end{tikzpicture}
	
\caption{\label{fig:sensitivity}Bidirectional bar charts, in which, for every parameter $p$, the value of $\sgn\left(\Upsilon^{\cR_0}_p\right)\left(18-r(p)\right)$ is represented by a bar which is labelled by the value of $\Upsilon^{\cR_0}_p$, in the disease-free case (top) and endemic case (bottom) specified in subsection \ref{subsec:sensitivity}. Bars to the left (right) of the ordinate axis are associated with parameters with which the basic reproduction number is monotonically decreasing (increasing). The longer the bar, the more significant the associated parameter. Consequently, the ordering of the parameters according to significance in each case is given by \eqref{eq:orderingDF} and \eqref{eq:orderingE}.}
\end{figure}
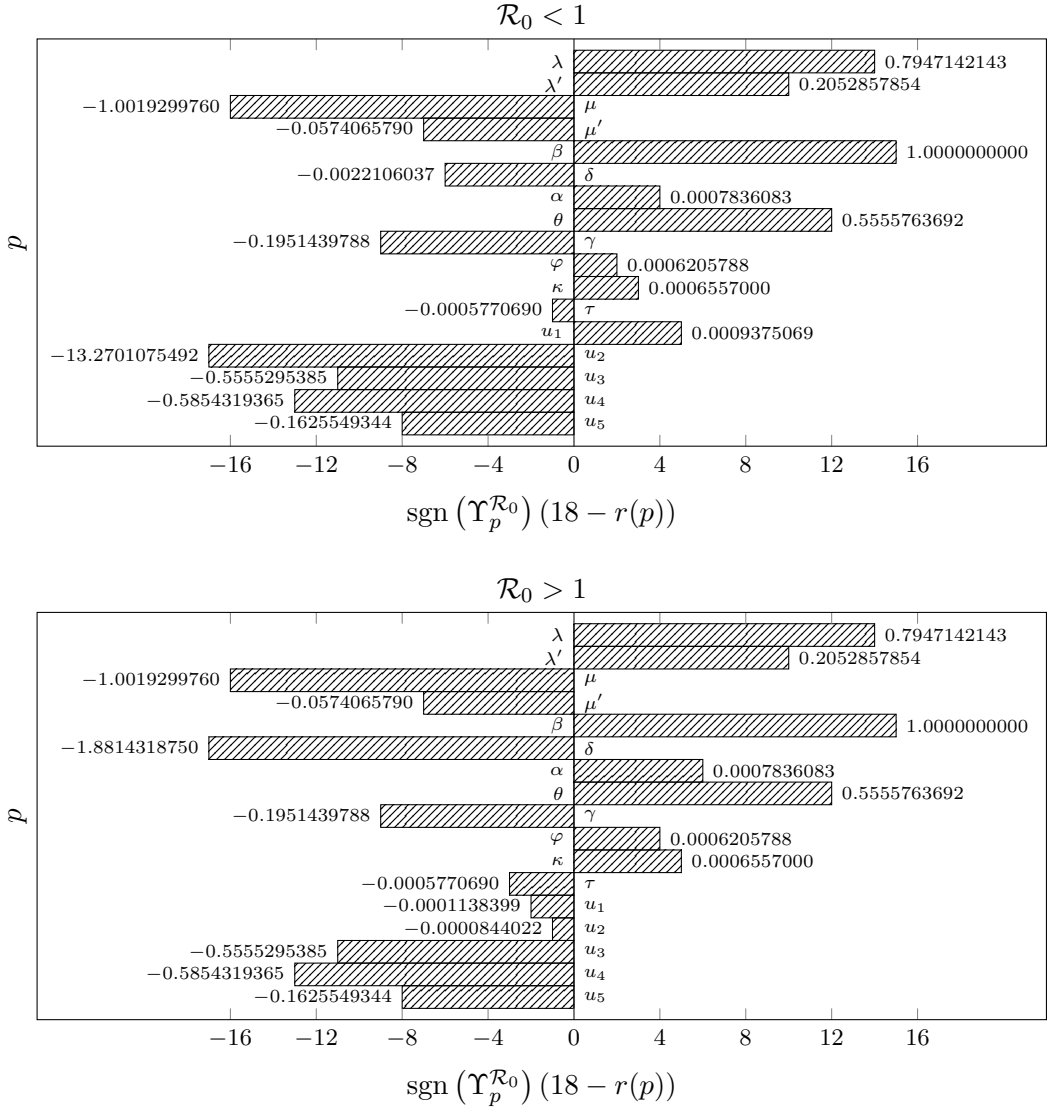

For a more effective visualisation, let us first define the \textit{significance rank} of a parameter $p$ to be $r(p)$, where $r:\mathcal{P}\to\{1,\ldots,17\}$ is the unique bijection for which
$$\left|\Upsilon^{\cR_0}_{r^{-1}(1)}\right|>\cdots> \left|\Upsilon^{\cR_0}_{r^{-1}(17)}\right|.$$
Thus, the values of $18-r(p)$ carry the same qualitative information as $\left|\Upsilon^{\cR_0}_p\right|$: the higher the value of $18-r(p)$, the more significant the parameter $p$. Therefore, instead of visualising the values of $\Upsilon^{\cR_0}_p=\sgn\left(\Upsilon^{\cR_0}_p\right)\left|\Upsilon^{\cR_0}_p\right|$, we visualise the values of $\sgn\left(\Upsilon^{\cR_0}_p\right)\left(18-r(p)\right)$, by bars which are labelled by the associated values of $\Upsilon^{\cR_0}_p$ to retain the quantitative information (Figure \ref{fig:sensitivity}). We can see that the ordering of parameters according to significance is
\begin{equation}\label{eq:orderingDF}
\left(r^{-1}(n)\right)_{n=1}^{17}=\left(u_2,\mu,\beta,\lambda,u_4,\theta,u_3,\lambda',\gamma,u_5,\mu',\delta,u_1,\alpha,\kappa,\varphi,\tau\right)
\end{equation}
in the disease-free case, and
\begin{equation}\label{eq:orderingE}
\left(r^{-1}(n)\right)_{n=1}^{17}=\left(\delta,\mu,\beta,\lambda,u_4,\theta,u_3,\lambda',\gamma,u_5,\mu',\alpha,\kappa,\varphi,\tau,u_1,u_2\right)
\end{equation}
in the endemic case. Let us now infer from these orderings the appropriate strategies of intervention in each case.


\begin{figure}
\includegraphics[scale=0.475]{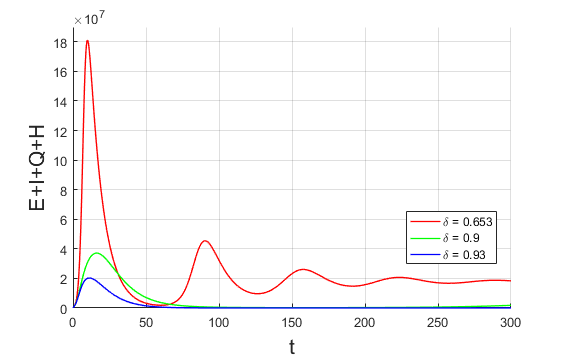} 
\caption{\label{fig:delta}Time-evolution of the number $E+I+Q+H$ of non-healthy individuals, for the parameter values shown in Table \ref{tab:parameters}, $u_1=0.4$, and $u_2=0.278$, with $\delta=0.653$ (red), $\delta=0.9$ (green), and $\delta=0.93$ (blue). Notice that, as $\delta$ is increased, both the maximum and the limit decrease drastically, showing the significance of $\delta$. For $\delta=0.93$, we observe convergence to the disease-free equilibrium.}
\end{figure}

\begin{figure}
\begin{tabular}{cc}(1)&(2)\\
\includegraphics[scale=0.475]{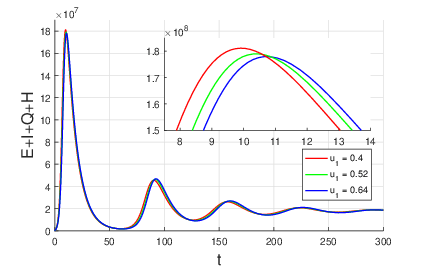} & \includegraphics[scale=0.475]{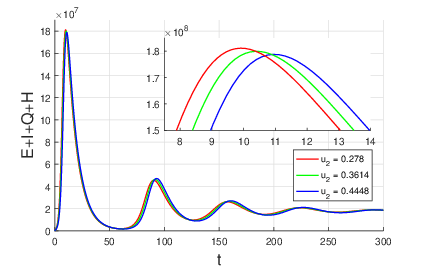}\\[0.25cm]
(3)&(4)\\
\includegraphics[scale=0.475]{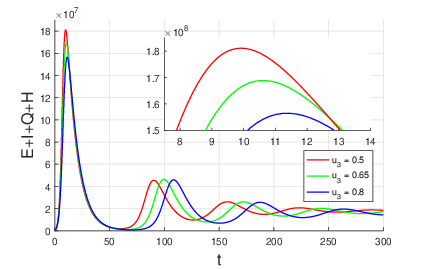} & \includegraphics[scale=0.475]{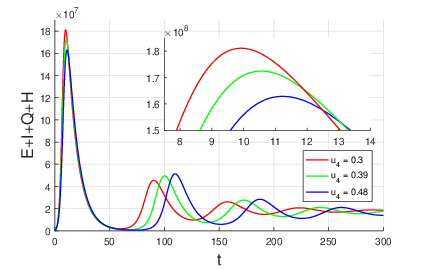}\\[0.25cm]
\multicolumn{2}{c}{(5)}\\
\multicolumn{2}{c}{
\includegraphics[scale=0.475]{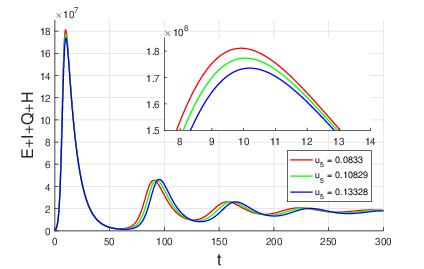}}
\end{tabular}
\caption{\label{fig:ui}Time-evolution of the number $E+I+Q+H$ of non-healthy individuals, for the parameter values shown in Table \ref{tab:parameters}, $\delta=0.653$, $u_1=0.4$, $u_2=0.278$, $u_3=0.5$, $u_4=0.3$, and $u_5=0.0833$ (red graphs on all panels), together with, on panel ($i$), the same in the cases of $u_i$ being increased by 30\% (green) and by 60\% (blue). The magnifications near the maxima reflect the previously obtained order of importance of the intervention parameters: $u_4$, $u_3$, $u_5$, $u_1$, $u_2$. However, comparing with Figure \ref{fig:delta}, we see that all these parameters are far less significant than $\delta$. Thus, the key for a successful eradication of the pandemic is not a high rate of implementation of any of the five intervention forms, but a high-efficacy vaccine.}
\end{figure}

\subsubsection{The disease-free case} In the disease-free case, effort must be made in order to maintain the low value of the basic reproduction number. When the number of daily new cases rises, as the ordering \eqref{eq:orderingDF} and Figure \ref{fig:sensitivity} suggest, a committed implementation of social restrictions should be sufficient. Indeed, remarkably, tightening social restrictions only by $10\%$ suppresses the basic reproduction number by $132.7010754919\%$. Other forms of intervention, if at all desired, are recommended in the following order: testings, tracings, treatments, and vaccinations. The latter is rather inessential, let alone when not supported by a high level of vaccine efficacy.

\subsubsection{The endemic case} In the endemic case, significant effort is necessary for a transition to a disease-free state. The ordering \eqref{eq:orderingE} strongly supports our main finding in subsection \ref{subsec:R0}: the vaccine efficacy being the parameter upon which the basic reproduction number depends most sensitively. Accordingly, we reiterate our primary finding in the previous subsection: that
\begin{enumerate}
\item raising the efficacy of vaccines
\end{enumerate}
must be given the highest priority. Indeed, increasing the efficacy from 65.3\%, firstly to 90\%, and subsequently to 93\%, as narrated in subsection \ref{subsec:R0}, results in significant drops of both the peak and limiting numbers of non-healthy individuals, the final value being sufficient for a transition to disease-free; see Figure \ref{fig:delta}.

With regards to the five forms of intervention, we recommend, in order of importance:
\begin{enumerate}[resume]
\item expanding and accelerating testings and tracings, so that infected and exposed individuals may be quarantined more immediately;
\item optimising treatments for infected individuals, by ensuring that health facilities and services (medications, hospital beds, medical practitioners, etc.)\ are in adequate availability;
\item administering vaccinations;
\item social restrictions, being the least important form of intervention, albeit, as previously remarked, may become significant if also applied to some degree to vaccinated individuals.
\end{enumerate}
There is however a large difference in significance between recommendations (1) and (2)--(5), i.e., the rates of all five forms of intervention are far less significant than the vaccine efficacy. This is apparent in Figure \ref{fig:ui}, where we can see that, assuming the original value of the vaccine efficacy, $65.3\%$, increasing any of the five intervention rates by 30\%, or even by 60\%, gives rise to barely any tangible impact.

\section{Conclusions and future research}\label{section:conclusions}

We have constructed an SVEIQHR-type mathematical model for the spread of COVID-19, which incorporates as parameters the rates of the five forms of intervention presently realised by the government of Indonesia: vaccinations, social restrictions, tracings, testings, and treatments. We have computed the model's basic reproduction number $\cR_0$, and show that the model possesses a unique disease-free equilibrium, which exists for all sets of parameter values and is stable (unstable) if $\cR_0<1$ ($\cR_0>1$), as well as a unique endemic equilibrium, which exists if $\cR_0>1$.

We have also analysed the model numerically, with the aim of determining strategies by which the five intervention forms should be realised in order to optimise their impact. The analysis results in two major conclusions. Firstly, in a disease-free state, social restrictions proved to be the best form of intervention in the case of a rise in the number of new cases. Secondly, in an endemic state, a transition to disease-free state without vaccinations requires a near-lockdown policy. Since the country's government has refused to impose such a policy \cite{GorbianoSutrisno,BorneoBulletin}, vaccinations are necessary. However, efforts should be focused not primarily on increasing the vaccination rate (or even the rate of any other form of intervention), but on the use of vaccines with a high efficacy.

Finally, our model is open to a number of modifications. One could incorporate a specified level of social restrictions for vaccinated individuals, and confirm whether, as a result, vaccination rate becomes more significant. Besides, the vaccinated compartment itself could be split into several compartments, in order to allow different assumptions on recipients of different vaccines and/or, in the case of multi-dose vaccines, recipients of different numbers of vaccine doses. Similarly, the quarantined compartment could be split into several compartments, in order to distinguish the isolated individuals (i.e., the separated infected individuals) from the quarantined individuals (i.e., the separated exposed individuals) ---which may further be split according to their vaccination histories--- so that different recovery rates and/or time-delays may be employed.


\end{document}